\documentclass{article}
\usepackage{amssymb,amsmath, euscript}

\input{epsf}

\newcounter{sec}

\newcounter{punct}[sec]

\def\punct{\refstepcounter{punct}{\arabic{sec}.\arabic{punct}.  }}

\def\COUNTERS{\addtocounter{sec}{1}
              \setcounter{punct}{0}
          \setcounter{equation}{0}
          \setcounter{theorem}{0}
               }

\newtheorem{theorem}{Theorem}[sec]
\newtheorem{proposition}[theorem]{Proposition}

\newtheorem{lemma}[theorem]{Lemma}

\newtheorem{corollary}[theorem]{Corollary}
\newtheorem{observation}[theorem]{Observation}

\newcounter{addendums}

 \newcounter{Apunct}

\begin{document}

\def\ov{\overline}

\def\wt{\widetilde}

\newcommand{\rk}{\mathop {\mathrm {rk}}\nolimits}
\newcommand{\Aut}{\mathop {\mathrm {Aut}}\nolimits}
\newcommand{\Out}{\mathop {\mathrm {Out}}\nolimits}
\renewcommand{\Re}{\mathop {\mathrm {Re}}\nolimits}

\def\Br{\mathrm {Br}}

\def\SL{\mathrm {SL}}
\def\SU{\mathrm {SU}}
\def\GL{\mathrm  {GL}}
\def\U{\mathrm  U}
\def\OO{\mathrm  O}
\def\Sp{\mathrm  {Sp}}
\def\SO{\mathrm  {SO}}
\def\SOS{\mathrm {SO}^*}
\def\Diff{\mathrm{Diff}}
\def\Vect{\mathfrak{Vect}}

\def\Ass{\mathrm{Asfr}}

\def\PGL{\mathrm  {PGL}}
\def\PU{\mathrm {PU}}

\def\PSL{\mathrm  {PSL}}

\def\Br{\mathrm{Br}}
\def\Fr{\mathrm{Fr}}

\def\Symp{\mathrm{Symp}}

\def\End{\mathrm{End}}
\def\Hom{\mathrm{Hom}}
\def\Mor{\mathrm{Mor}}
\def\Aut{\mathrm{Aut}}

\def\PB{\mathrm{PB}}

\def\cA{\mathcal A}
\def\cB{\mathcal B}
\def\cC{\mathcal C}
\def\cD{\mathcal D}
\def\cE{\mathcal E}
\def\cF{\mathcal F}
\def\cG{\mathcal G}
\def\cH{\mathcal H}
\def\cJ{\mathcal J}
\def\cI{\mathcal I}
\def\cK{\mathcal K}
\def\cL{\mathcal L}
\def\cM{\mathcal M}
\def\cN{\mathcal N}
\def\cO{\mathcal O}
\def\cP{\mathcal P}
\def\cQ{\mathcal Q}
\def\cR{\mathcal R}
\def\cS{\mathcal S}
\def\cT{\mathcal T}
\def\cU{\mathcal U}
\def\cV{\mathcal V}
\def\cW{\mathcal W}
\def\cX{\mathcal X}
\def\cY{\mathcal Y}
\def\cZ{\mathcal Z}

\def\sm{\smallskip}


\def\0{{\ov 0}}
\def\1{{\ov 1}}


\def\frA{\mathfrak A}
\def\frB{\mathfrak B}
\def\frC{\mathfrak C}
\def\frD{\mathfrak D}
\def\frE{\mathfrak E}
\def\frF{\mathfrak F}
\def\frG{\mathfrak G}
\def\frH{\mathfrak H}
\def\frJ{\mathfrak J}
\def\frK{\mathfrak K}
\def\frL{\mathfrak L}
\def\frM{\mathfrak M}
\def\frN{\mathfrak N}
\def\frO{\mathfrak O}
\def\frP{\mathfrak P}
\def\frQ{\mathfrak Q}
\def\frR{\mathfrak R}
\def\frS{\mathfrak S}
\def\frT{\mathfrak T}
\def\frU{\mathfrak U}
\def\frV{\mathfrak V}
\def\frW{\mathfrak W}
\def\frX{\mathfrak X}
\def\frY{\mathfrak Y}
\def\frZ{\mathfrak Z}

\def\fra{\mathfrak a}
\def\frb{\mathfrak b}
\def\frc{\mathfrak c}
\def\frd{\mathfrak d}
\def\fre{\mathfrak e}
\def\frf{\mathfrak f}
\def\frg{\mathfrak g}
\def\frh{\mathfrak h}
\def\fri{\mathfrak i}
\def\frj{\mathfrak j}
\def\frk{\mathfrak k}
\def\frl{\mathfrak l}
\def\frm{\mathfrak m}
\def\frn{\mathfrak n}
\def\fro{\mathfrak o}
\def\frp{\mathfrak p}
\def\frq{\mathfrak q}
\def\frr{\mathfrak r}
\def\frs{\mathfrak s}
\def\frt{\mathfrak t}
\def\fru{\mathfrak u}
\def\frv{\mathfrak v}
\def\frw{\mathfrak w}
\def\frx{\mathfrak x}
\def\fry{\mathfrak y}
\def\frz{\mathfrak z}

\def\frsp{\mathfrak{sp}}

\def\frfr{\frf\frr}
\def\frbr{\frb\frr}


\def\bfa{\mathbf a}
\def\bfb{\mathbf b}
\def\bfc{\mathbf c}
\def\bfd{\mathbf d}
\def\bfe{\mathbf e}
\def\bff{\mathbf f}
\def\bfg{\mathbf g}
\def\bfh{\mathbf h}
\def\bfi{\mathbf i}
\def\bfj{\mathbf j}
\def\bfk{\mathbf k}
\def\bfl{\mathbf l}
\def\bfm{\mathbf m}
\def\bfn{\mathbf n}
\def\bfo{\mathbf o}
\def\bfp{\mathbf p}
\def\bfq{\mathbf q}
\def\bfr{\mathbf r}
\def\bfs{\mathbf s}
\def\bft{\mathbf t}
\def\bfu{\mathbf u}
\def\bfv{\mathbf v}
\def\bfw{\mathbf w}
\def\bfx{\mathbf x}
\def\bfy{\mathbf y}
\def\bfz{\mathbf z}

\def\bfA{\mathbf A}
\def\bfB{\mathbf B}
\def\bfC{\mathbf C}
\def\bfD{\mathbf D}
\def\bfE{\mathbf E}
\def\bfF{\mathbf F}
\def\bfG{\mathbf G}
\def\bfH{\mathbf H}
\def\bfI{\mathbf I}
\def\bfJ{\mathbf J}
\def\bfK{\mathbf K}
\def\bfL{\mathbf L}
\def\bfM{\mathbf M}
\def\bfN{\mathbf N}
\def\bfO{\mathbf O}
\def\bfP{\mathbf P}
\def\bfQ{\mathbf Q}
\def\bfR{\mathbf R}
\def\bfS{\mathbf S}
\def\bfT{\mathbf T}
\def\bfU{\mathbf U}
\def\bfV{\mathbf V}
\def\bfW{\mathbf W}
\def\bfX{\mathbf X}
\def\bfY{\mathbf Y}
\def\bfZ{\mathbf Z}

\def\bfBr{\bfB\bfr}
\def\bfFr{\bfF\bfr}

\def\G{\mathbf G}

\def\bfw{\mathbf w}

\def\R {{\mathbb R }}
 \def\C {{\mathbb C }}
  \def\Z{{\mathbb Z}}
  \def\H{{\mathbb H}}
\def\K{{\mathbb K}}
\def\N{{\mathbb N}}
\def\Q{{\mathbb Q}}
\def\A{{\mathbb A}}

\def\T{\mathbb T}
\def\P{\mathbb P}

\def\bbA{\mathbb A}
\def\bbB{\mathbb B}
\def\bbD{\mathbb D}
\def\bbE{\mathbb E}
\def\bbF{\mathbb F}
\def\bbG{\mathbb G}
\def\bbI{\mathbb I}
\def\bbJ{\mathbb J}
\def\bbL{\mathbb L}
\def\bbM{\mathbb M}
\def\bbN{\mathbb N}
\def\bbO{\mathbb O}
\def\bbP{\mathbb P}
\def\bbQ{\mathbb Q}
\def\bbS{\mathbb S}
\def\bbT{\mathbb T}
\def\bbU{\mathbb U}
\def\bbV{\mathbb V}
\def\bbW{\mathbb W}
\def\bbX{\mathbb X}
\def\bbY{\mathbb Y}

\def\kappa{\varkappa}
\def\epsilon{\varepsilon}
\def\phi{\varphi}
\def\le{\leqslant}
\def\ge{\geqslant}
\def\la{\langle}
\def\ra{\rangle}

\def\tr{\mathrm{tr}\,}

\def\E{\cE}

\begin{center}

{\bf \Large On topologies on Malcev completions of braid groups}

\bigskip

{\sc\large Yury Neretin%
\footnote{Supported by grant FWF P22122.} }

\end{center}

\begin{flushright}
 To the memory of Vladimir Igorevich Arnold
\end{flushright}

{\small We discuss groups corresponding to Kohno Lie algebra of
infinitesimal braids and  actions of such groups.
 We construct
homomorphisms of Lie braid groups
 to the group of symplectomorphisms
of the space of point configurations in $\R^3$ and to groups of
symplectomorphisms of coadjoint orbits of $\SU(n)$.}

\section{Introduction. Malcev completion of the group of pure braids}

\COUNTERS

 Artin braid groups 
are well-known objects, which appear in natural way in many branches of mathematics
and mathematical physics. 
 According an abstract Malcev construction the group $\Br_n$ of pure braids admits
 a canonical embedding to a certain infinite-dimensional nilpotent Lie group $\bfBr_n$.
  The Lie algebra $\frbr_n$ of $\bfBr_n$  described by
 Kohno also was a topic 
 of numerous investigations.
 
 The first (informal) purpose of this paper is to explain that the infinite-dimensional Lie group $\bfBr_n$ is an interesting and relatively hand-on object.  The formal purpose is to 
 prove several statements about these groups.
 
 \sm

{\bf\punct Group of pure braids.} Denote by $\Br_n$ the group of
pure braids of $n$ strings (see e.g, \cite{Art}, \cite{Mag},
\cite{Bir}). It is a discrete group with generators
 $A_{rs}$, where $1\le r<s\le n$, and relations
 \begin{align*}
&\{A_{rs},A_{ik}\}=1,\quad &\text{if $s<i$ or $k<r$}
\\
&\{A_{ks},A_{ik}\}=\{A_{is}^{-1},A_{ik}\},\quad &\text{if $i<k<s$}
\\
&\{A_{rk},A_{ik}\}=\{A_{ik}^{-1}, A_{ir}^{-1}\}
    ,\quad &\text{if $i<r<k$}
    \\
& \{A_{rs},A_{ik}\}=
  \bigl\{\{ A_{is}^{-1},A_{ir}^{-1}\},A_{ik}\bigl\}
  ,\quad &\text{if $i<r<k<s$}
 \end{align*}
Here $\{a,b\}:=aba^{-1}b^{-1}$ denotes a {\it commutator}. For geometric interpretation
of the generators, see \cite{Bir}, Fig.4.

There is a natural homomorphism $\Br_{n}\to \Br_{n-1}$, we forget
the first string. It is easy to verify, that
 the kernel is the {\it free
group} $F_{n-1}$ of $(n-1)$ generators. Therefore,
 $\Br_n$ is a
semi-direct product
$$
\Br_n=\Br_{n-1}\ltimes F_{n-1}
.$$
Repeating the same argument, we get that
$\Br_n$ is a product of its subgroups
$$
\Br_n=F_1\times F_2\times \dots\times F_{n-1} \qquad
 \text{as a set,}
$$
i.e., each $g\in\Br_n$ admits a unique representation
$g=h_{n-1}\dots h_2 h_1$, where $h_j \in F_j$ (we prefer the inverse order).
 More precisely,
$$
\Br_n=F_1\ltimes \Bigl( F_2\ltimes\Bigl(F_3\ltimes\dots
\ltimes \Bigl( F_{n-2}\ltimes F_{n-1}\Bigr)\Bigr)\Bigr)
$$
as a group and subgroups $F_k\ltimes( \dots\ltimes F_{n-1})$ are normal.

\smallskip


{\bf\punct Abstract Malcev construction (pro-unipotent
completion). Step 1.\label{ss:maltsev}}
Let $\Gamma$ be a group, $A$, $B$ two subgroups. Denote by
$\{A,B\}$ their commutant, i.e., the subgroup generated by all
commutators $\{a,b\}$, where $a\in A$, $b\in B$.
 For a group $\Gamma$ denote
$\Gamma_1:=\{\Gamma,\Gamma \}$, and $\Gamma_{j+1}:=\{\Gamma_j,\Gamma\}$.
A discrete  group $\Gamma$ is  {\it nilpotent} if some
$\Gamma_j=1$ and
 {\it residually nilpotent} if $\cap_j \Gamma_j=1$. The group
$\Br_n$ is residually nilpotent.

Malcev (see \cite{Mal1} and comprehensive  expositions in \cite{Rag}, Chapter II, and
 \cite{KM}, Section 17)
 proved that any nilpotent discrete group
$\Gamma$ without torsion admits a canonical embedding
$\Gamma\to G$ to a
 simply connected nilpotent Lie group $G$ as
a uniform lattice%
\footnote{i.e, $\Gamma$ is a discrete subgroup and
 the homogeneous space $G/\Gamma$ is
compact}.
He constructs 'coordinates' on $\Gamma$
in the following way. Choose a basis in each free Abelian
group $\Gamma_{j-1}/\Gamma_j$, take representatives
$a_1$, $a_2$, \dots,  $a_N$ of
this basis elements in $\Gamma$, let an initial segment
of this sequence  corresponds
to $\Gamma/\Gamma_1$, a next segment corresponds to
$\Gamma_1/\Gamma_2$ etc., Any element of $\Gamma$ admits a
unique representation in the form
$a_1^{q_1} a_2^{q_2}\dots a_N^{q_N}$. Therefore, we
can identify sets $\Gamma$ and $\Z^N$. Malcev proves
that formulas
for multiplication in 'coordinates' $(q_1,\dots,q_N)$
are polynomial. We  take real $q_j$ and get a nilpotent
Lie group $G$, where a multiplication
is given by the same formulas.

This construction is functorial ('rigid') in the following senses. 

\sm

$1^\circ$. {\it  Let $\Gamma$, $\Delta$ be nilpotent groups without torsion,
and $\rho:\Gamma\to \Delta$ be a homomorphism. Then $\rho$ admits a unique
extension to the homomorphism $G\to D$ of the corresponding completions.}

\sm

$2^\circ$. {\it  Let $\Gamma$ be a nilpotent group without torsion, $G$ be its completion.
Let $H$ be a simply connected nilpotent Lie group. Then any homomorphism $\Gamma\to H$ admits a unique
continuous extension to a homomorphism $G\to H$}
(see \cite{Rag}, Theorem 2.11).

\sm

$3^\circ$. {\it Let $\rho:\Gamma\to \GL(n,\C)$ be a homomorphism such that all matrices
$\rho(\gamma)$ be unipotent%
\footnote{A matrix is unipotent if all its eigenvalues equal 1.}.
Then $\rho$ admits a unique continuous extension to a homomorphism $G\to\GL(n,\C)$}.

\sm

{\sc Remark.}
Evidently, the $2^\circ$ implies $1^\circ$. On the other hand $2^\circ$ implies $3^\circ$. Indeed, consider the Zariski  closure $\ov{\rho(\Gamma)}$ of $\rho(\Gamma)$. First,  for any unipotent matrix $A\ne 1$ the Zariski closure of the set $\{A^n\}$ is the unipotent one-parametric subgroup 
$$
A^s=\exp\Bigl\{s\sum_{j\ge 0} \frac{(-1)^{j-1}}j (A-1)^j\Bigr\}, \qquad s\in \C
.$$
Therefore the subgroup $\ov{\rho(\Gamma)}$ is connected. On the other hand $\ov{\rho(\Gamma)}$
consists from unipotent matrices. By the Lie theorem  (see. e.g., \cite{Dix}, 1.3.12)
 it admits an invariant flag. Therefore we get a
homomorphism from $\Gamma$ to upper triangular matrices and can refer to $2^\circ$.

\smallskip

{\sc Remark.} There is another way of construction of the completion proposed by Jennings
\cite{Jen}, see also \cite{Qui}. Malcev \cite{Mal2}, \cite{Mal3} also constructed the group
of rational points of $G$.

\smallskip

{\bf\punct Step 2.}
Similarly, any residually nilpotent group
without torsion  admits a canonical embedding to a certain  projective limit
of nilpotent Lie groups. Indeed, consider groups $\Gamma/\Gamma_j$ and their completions
$G^{(j)}$. We have a chain of homomorphisms 
\begin{equation}
\dots
\leftarrow \Gamma/\Gamma_j \leftarrow \Gamma/\Gamma_{j+1}\leftarrow\dots
\label{eq:chain-1}
\end{equation}
and therefore a chain of homomorphisms
$$
\dots
\leftarrow G^{(j)} \leftarrow G^{(j+1)} \leftarrow\dots
$$
We denote by $\ov\bfG$ its projective limit. By the construction
$\Gamma$ is a subgroup of $\ov\bfG$.

On the other hand, the inverse limit $\ov{\mathbf{\Gamma}}\subset \ov\bfG$ 
of the chain (\ref{eq:chain-1}) 
is a continual topological group. It contains $\Gamma$ as a dense  subgroup.

\smallskip

{\bf\punct Kohno algebra.%
\label{ss:kohno}}
The braid group  $\Br_n$ is residually nilpotent, therefore it has
a completion $\bfBr_n$ and the corresponding Lie algebra $\frbr_n$.

Kohno \cite{Koh} obtained an explicit description of
the Lie algebra $\frbr_n$.
It is generated by elements $r_{ij}$, where
$0\le i,j\le n$  and $r_{ij}=r_{ji}$, the relations are
\begin{align}
&[r_{ij},r_{kl}]=0,\qquad\text{if $i$, $j$, $k$, $l$ are pairwise
distinct,}
\label{eq:koh1}
\\
&[r_{ij},r_{ik}+r_{jk}]=0.
\label{eq:koh2}
\end{align}

 The Campbell--Hausdorff formula (see an introduction in Serre
\cite{Ser} and a treatise in Reutenauer \cite{Reu})
 produces a group structure on $\frbr_n$,
$$
x\cdot y:=\ln( e^x e^y).
$$
This is the corresponding {\it Lie braid group $\ov\bfBr_n$}.

The algebra $\frbr_n$ is $\Z_+$-graded, the gradation is defined 
 by the condition $\deg r_{ij}=1$.
Kohno obtained the Poincare series  for
dimensions of homogeneous subspaces of the universal enveloping
algebra $\U(\frbr_n)$ of $\frbr_n$,
\begin{equation}
\sum_{k\ge 0}^\infty t^k\cdot \dim \U(\frbr_n)^{[k]}=
\prod_{j=1}^{n-1}(1-jt)^{-1}
.\label{eq:koh-poi}
\end{equation}
This easily implies that
$$
 \dim \frbr_n^{[k]}
=\frac 1k \sum_{d|k} \mu(d)\Bigl( \sum_{i=1}^{n-1} i^{k/d}\Bigr)
,$$
where the summation is given over all divisors $d$ of $n$ and
$\mu(d)$ is the M\"obius function%
\footnote{if $d=p_1\dots p_l$ is
a product of pairwise distinct primes, then $\mu(d)=(-1)^l$,
otherwise $\mu(d)=0$.} (proof as in \cite{Ser}, Section IV.4).
In particular, $\dim \frbr_n^{[k]}$ has exponential growth
 as $k\to\infty$.

Also, Kohno noted (see \cite{Koh2}, \cite{Koh3},
 and below \ref{ss:KZ})
 that any finite-dimensional representation
of the Lie algebra $\frbr_n$ produces a representation
 of the braid group
in the same space.

\smallskip

{\bf\punct Knizhnik--Zamolodchikov construction.\label{ss:KZ1}} Apparently 
the most important origin
of representations of $\frbr_n$ is the following construction.
Let $\frg$ be a semisimple Lie algebra.
 Consider the enveloping algebra
$\U(\frg\oplus\dots\oplus\frg)$ of direct sum of $n$ copies of $\frg$
Let $e_\alpha$ be an orthormal basis
in $\frg$. Let $e_\alpha^j$ be the same basis in $j$-th copy of $\frg$.
Then mixed Casimirs
$$
\Delta_{ij}=\sum_\alpha e^i_\alpha e^i_\beta
$$
satisfy the relations (\ref{eq:koh1})--(\ref{eq:koh2}).

Therefore $\frbr_n$ acts in any tensor product $V_1\otimes\dots\otimes V_n$
of finite-dimensional representations of $\frg$.

Many  finite-dimensional representations of $\frbr_n$ are known, we refer to 
\cite{EFK}, \cite{Koh2}, \cite{Tol}, \cite{Mil}.

\smallskip


{\bf\punct Knizhnik--Zamolodchikov connection.\label{ss:KZ2}} See \cite{Koh2}, \cite{Koh3}.
 Now consider a collection of 
linear operators $r_{ij}$ in a finite-dimensional space 
$V$ satisfying the commutation relations  (\ref{eq:koh1})--(\ref{eq:koh2}).
Consider the space $\Omega$ obtained from $\C^n$ by the removing of all diagonals
$z_k=z_l$. Its fundamental group $\pi_1(\Omega)$ is $\Br_n$.
 Consider the connection on $\Omega$ given by the formula
\begin{equation}
\sum_{k>l} r_{kl}\cdot d\ln(z_i-z_j)
.
\label{eq:KZ}
\end{equation}
A straightforward calculation shows that the connection is flat.
The monodromy of the connection determines a linear representation of the group $\pi_1(\Omega)\simeq\Br_n$
in $V$.

\sm


{\bf\punct Klyachko's spatial polygons.\label{ss:klyachko}}
 Consider the space
$\R^3$  and the Poisson bracket on the space of functions determined by
$$
\{x,y\}=z,\qquad \{y,z\}=x,\qquad \{z,x\}=y
$$
or, more precisely,
$$
\{f,g\}= z\,
\Bigl(\frac{\partial f}{\partial x}\cdot \frac{\partial g}{\partial y}
-
\frac{\partial f}{\partial y}\cdot \frac{\partial g}{\partial x}\Bigr)
+
x\, \Bigl(
\frac{\partial f}{\partial y}\cdot \frac{\partial g}{\partial z}
-
\frac{\partial f}{\partial z}\cdot \frac{\partial g}{\partial y}\Bigr)
+
y\,\Bigl(\frac{\partial f}{\partial z}\cdot \frac{\partial g}{\partial x}
-
\frac{\partial f}{\partial x}\cdot \frac{\partial g}{\partial z}\Bigr)
.$$
Consider the space $(\R^3)^n$ whose points are collections
of 3-vectors $(r_1,\dots,r_n)$. Symplectic leafs of our Poisson structure are products of spheres
$|r_j|=a_j$.

 The following  collection of Hamiltonians
$$
\Delta_{ij}:=\la r_i,r_j\ra =x_i x_j+ y_iy_j+z_iz_j,
\qquad i\ne j,
$$
 satisfies the Kohno relations
(\ref{eq:koh1})--(\ref{eq:koh2}).

In fact, Klyachko \cite{ Klya} considered the symplectic structure on
the space of point configurations $(r_1,\dots,r_n)$ such that
$$
\sum r_j=0,\qquad |r_j|=a_j
$$
up to rotations. This is the space of spatial polygons with given lengths of
sides (this object has an unexpectedly rich geometry).

\smallskip


{\bf\punct Extension of the Klyachko construction.}
More generally, consider a semisimple Lie algebra $\frg$.
 Consider the standard Poisson bracket
on $\frg^*$. If  $u_\alpha\in \frg$ is a basis, we write
$$
\{f,g\}:= \sum_{i,j} [u_i,u_j]\frac{\partial f}{\partial u_i}
\frac{\partial g}{\partial u_j}
.$$
Next, consider the product $\frg^*\times\dots\times\frg^*$, we denote its points
as $(X_1,\dots,X_n)$.
Then the collection of Hamiltonians
$$
\Delta_{ij}:=\tr X_i X_j
$$
satisfies the Kohno relations (\ref{eq:koh1})-(\ref{eq:koh2}).

In Klyacko's example, $\frg=\mathfrak{so}_3$.

Evidently, this construction is a classical limit of the Knizh\-nik--Za\-mo\-lod\-chi\-kov construction.
 The next construction
is a classical limit of a construction discussed in  \cite{Mil}.

\smallskip

{\bf\punct Action of the Kohno algebra on coadjoint orbits.%
\label{ss:laredo}}
Consider the Lie algebra $\mathfrak{gl}(n)$ of all $n\times n$ matrices
 and its dual  $\mathfrak{gl}(n)^*$. Consider the
Poisson structure on the space of functions on $\mathfrak{gl}(n)^*$
(we denote elements of this space as matrices $\{x_{ij}\}$).
 Symplectic leafs are coadjoint orbits.
Then the collection of Hamiltonians
$$
\Delta_{ij}=2x_{ij}x_{ji}
$$
form a representation of the Kohno algebra.

{\sc Proof.} Elements $x_{ij}$ with $i$, $j\le 3$ form the Lie algebra
$\mathfrak{gl}(3)$. The expression
$$
\tr X^2=\sum_{i\le 3} X_{ii}^2+ 2\sum_{i<j\le 3} X_{ij}X_{ji}=\sum_{i\le 3} X_{ii}^2+
\Delta_{12} +\Delta_{13}+\Delta_{23}
$$
commutes with all $x_{ij}$ with $i$, $j\le 3$. Therefore it commutes with $\Delta_{12}$.
On the other hand $\Delta_{12}$ commutes with $X_{11}$, $X_{22}$, $X_{33}$. Hence 
$$
\{\Delta_{12},  \sum_{i\le 3} X_{ii}^2\}=0.
$$
Therefore
$$
\{\Delta_{12}, \Delta_{12} +\Delta_{13}+\Delta_{23}\}=0
.
$$

\smallskip

{\bf\punct An example.} Knizhnik--Zamoldchikov construction
produces numerous examples of actions of $\frbr_n$ by
second-order partial differential operators. For instance consider
the product of $n$ copies of $\R^2$ with coordinates
$(x_j,y_j)$. The operators
\begin{equation}
\Delta_{kl}=(x_ky_l-x_ly_k)\Bigl(
\frac{\partial^2}{\partial x_k\partial y_l}-
\frac{\partial^2}{\partial x_l\partial y_k}
\Bigr)
\label{eq:laplace}
\end{equation}
satisfy the same commutation relations (\ref{eq:koh1})--(\ref{eq:koh2}).

\smallskip

{\bf \punct Some relatives of the braid Lie algebra.\label{ss:misc}}
 The Kohno algebra is a representative
of a big zoo of natural infinite-dimensional nilpotent Lie algebras. 

\sm

a) {\it Generalized Kohno algebras}, see \cite{Mil}.
 We consider a reduced root system $R$. For each
positive root $\alpha$ (a notation: $\alpha\in R^+$),
 we assign a generator $r_\alpha$. For each subsystem $S\subset R$ of rank 2
 and $\alpha\in S^+$
 we assign a relation 
 $$
 [r_\alpha, \sum_{\gamma\in S} r_\gamma]=0
 .
 $$

b) {\it Groups related to Yang-Baxter equation,}
 see \cite{BER}, \cite{Lee}. Consider the '{\it quasitriangular} group' $\mathrm{QTr}_n$, it has generators
$R_{ij}$, where $1\le i,j\le n$, $i\ne j$, and relations
$$
 R_{ij} R_{ik} R_{jk} = R_{jk} R_{ik} R_{ij} .
$$
This is a residually nilpotent Lie group, its Lie algebra $\mathfrak{qtr}_n$
 is generated by
elements $r_{ij}$, where $1\le i,j\le n$, $i\ne j$, the relations are
$$
[r_{ij} , r_{ik} ] + [r_{ij} , r_{jk} ] + [r_{ik} , r_{jk} ] = 0
.
$$

The {\it triangular group} $\mathrm{Tr}_n$ is the group with the same generators and the additional relation
$R_{ij}=R_{ji}^{-1}$. For the corresponding Lie algebra $\mathfrak{qtr}_n$
 we have additional relation
$r_{ij}=-r_{ji}$.

By the definition, the algebra $\mathfrak{tr}_n$ is a quotient of $\mathfrak{qtr}_n$.
On the other hand the subalgebra in $\mathfrak{qtr}_n$ generated by all $r_{ij}+r_{ji}$
is the Kohno algebra $\frbr_n$.

\sm

c) {\it  Torelli groups.} We consider 2-dimensional compact surface $M$ of genus $g$ with $k$
distinguished points. Consider the group $Diff(M)$ of diffeomorphisms of $M$, fixing
the distinguished  points.
Denote by $Diiff_0(M)$ be the connected component of unit.  Consider the mapping class group
$Diff(M)/Diff_0(M)$. It acts on the first homology group $H_1(M,\Z)$ of $M$, the group
$H_1(M,\Z)$ is equipped with a skew symmetric bilinear form (the intersection index).
Therefore we get a homomorphism $Diff(M)/Diff_0(M)\to \Sp(2n,\Z)$. Kernels of such homomorphisms are  residually nilpotent, their Lie algebras were described by Hain,
\cite{Hai} (he also described  completions of the whole groups $Diff(M)/Diff_0(M)$).

The braid group $\Br_n$ corresponds to the case $M=\R^2$.


{\bf \punct Purposes of the paper.}
As we have seen (\ref{ss:maltsev}), (\ref{ss:KZ2})
 there are correspondences between finite-dimensional
 representations of
the braid group $\Br_n$ and representations of the Kohno algebra $\frbr_n$.
However, representations of $\frbr_n$ do not produce representations of
the Lie braid group $\ov\bfBr_n$ (since formal series  diverge).

There is a well-known way to avoid convergence problems by passing
to spaces over algebra of formal series. For instance, in the Klyachko case,
we add a formal variable $\hbar$ and take new generators $\wt\Delta_{kl}:=\hbar\Delta_{kl}$.
Then any formal series in $\wt\Delta_{kl}$ determines a well-defined operator
in the space  $C^\infty((\R^3)^n)\otimes \C[\hbar]$. 

However Klyachko's Hamiltonians $\Delta_{kl}$
and their linear combinations generate nice Hamiltonian flows,
 there arises a question
about the group generated by all such flows.

\smallskip

In this paper, we define two topological versions of Lie braid groups 
 $\bfBr^\circ_n$, $\bfBr^!_n$, which are dense subgroups
of $\ov\bfBr_n$. They satisfy the following properties

--- any representation of  the Lie algebra $\frbr_n$  in a finite-dimensional
linear space can be integrated to
a representation of the group $\bfBr^\circ_n$ (see Section 4).

--- any action of the Lie algebra $\frbr_n$ on a compact real analytic
manifold can be integrated to an action of
 $\bfBr^!_n$ (Sections 5--6).
 
 The most of proofs are simple, an exception is Theorem \ref{th}.

I also try to explain that the Lie braid groups are
hand-on (or at least semi-hand-on) objects (Sections 2--3 and Theorem \ref{th:TTT}).
A rich theory of free Lie algebras (and 'free Lie groups'), see, e.g., the book
of Reutenauer \cite{Reu} can be
applied as tool in our situation (but the paper contains only few steps in this direction).

\smallskip

In fact, we consider infinite-dimensional nilpotent Lie groups in a wider generality,
this generality includes to considerations all the 'relatives' of the Kohno algebra mentioned
in \ref{ss:misc}. Also we discuss topological completions of  enveloping  algebras
in spirit of Rashevskij \cite{Rash}, see Propositions \ref{l:Ucirc}, \ref{l:!}.

\smallskip

We can not formulate any positive statement related to operators
(\ref{eq:laplace}) or for Knizhnik--Zamolodchikov actions of $\frbr_n$ in tensor products of unitary
representations of semisimple Lie groups (an exception is a semi-obvious Proposition 
\ref{pr:highest}).


{\bf Acknowledgements.} I am grateful to I.~A.~Dynnikov,
R.~S.~Ismagilov,
S.~M.~Khoroshkin, A.~M.~Levin, P.~Michor,
A.~A.~Rosly, L.~G.~Rybnikov for discussions of this topics.

\section{Groups $\ov \G$}

\COUNTERS

This section contains preliminary formalities.

\smallskip

{\bf\punct Free associative algebra.} Denote by $\Ass_n=\Ass[\omega]$
the free associative algebra with generators
$\omega_1$,\dots, $\omega_n$. The monomials
$w=\omega_{i_1}\dots\omega_{i_m}$ form a basis in $\Ass[\omega]$.
By $\ov\Ass_n[\omega]$ we denote the corresponding  algebra of 
formal series.

Denote by $\frfr_n=\frfr[\omega]$ the free Lie algebra with
generators $\omega_j$. Then (see, e.g, \cite{Ser}, \cite{Reu})
$\Ass_n$ is the universal enveloping algebra of $\frfr_n$,
$$
\Ass[\omega]\simeq \U(\frfr[\omega]).
$$


{\bf\punct Completed enveloping algebras.}
Let $\frg$ be a graded real Lie algebra
$$
\frg=\bigoplus\limits_{k>0} \frg^{[k]}, \qquad
[\frg^{[k]}, \frg^{[m]}]\subset \frg^{[k+m]}
,$$
let $\dim \frg^{[k]}<\infty$, assume that $\frg^{[1]}$ generates
$\frg$. We also fix a basis $\omega_i\in \frg^{[1]}$
(our main objects, groups $\ov \G$, $\G^\circ$, and $\G^!$,
 do not depend on a choice of basis).

\sm

{\sc Examples.}
1. Our basic example is the Kohno Lie algebra $\frbr_n$.

2) Another
example under discussion is the free Lie algebra  $\frfr_m$ with
$m$ generators

3) Arbitrary Lie algebra from \ref{ss:misc}.\hfill $\square$

\sm

Denote by $\ov\frg$ the Lie algebra whose elements
 are formal series
 $$
x=x^{[1]}+x^{[2]}+\dots,\qquad\text{where $x^{[j]}\in\frg^{[j]}$}
. $$

Denote by $\U(\frg)$ the universal enveloping algebra of $\frg$
equipped with the natural gradation arising from $\frg$.
We expand an element $z\in \U(\frg)$ as a sum of homogeneous
elements
\begin{equation}
z=z^{[0]}+z^{[1]}+z^{[2]}+z^{[3]}+\dots,\qquad
\qquad \text {where $z^{[j]}\in \U(\frg)^{[j]}$}
.
\label{z}
\end{equation}
Such sums are finite.
Denote by $\ov \U (\frg)$ the space
 of infinite formal  series
(\ref{z}). We equip the space $\ov\U(\frg)$ with the topology
of direct product of finite-dimensional linear spaces,
 $\ov\U(\frg)=\prod_j \U(\frg)^{[j]}$.

Denote by $\ov \U_+(\frg)$ the space of all $z$ such that
$ z^{[0]}=0$.

\begin{lemma}
{\rm a)} $\ov \U(\frg)$ is an algebra.

\smallskip

{\rm b)} For any formal series $h(\xi)=\sum_{j\ge 0} c_j \xi^j$,
where $c_j\in \C$, and any $z\in \ov \U_+(\frg)$,
the formal series $h(z)\in \ov\U(\frg)$ is well defined.

\smallskip

{\rm c)} For any $z\in \ov \U_+(\frg)$ the exponent $\exp(z)$ and the
logarithm $\ln (1+z)$ are well defined.
\end{lemma}

This is obvious.

\smallskip


{\bf\punct Group $\ov \G$.}

\begin{proposition}
The set $\ov\G$ of all elements $\exp(x)\in\ov \U(\frg)$,
 where $x$ ranges in $\frg$, is a group with respect to multiplication.
\end{proposition}

This follows from the Campbell--Hausdorff formula, see, e.g.,
\cite{Ser} or \cite{Reu}.

\smallskip


{\bf\punct Characterization of $\ov \G$.}
Consider the tensor product $\ov \U(\frg)\otimes \ov \U(\frg)$.
Consider the homomorphism ({\it co-product})
$\delta:\ov \U(\frg)\to\ov \U(\frg)\otimes \ov \U(\frg)$ by
$$
\delta(x)=x\otimes 1+1\otimes x, \qquad\text{for $x\in\frg$}
$$

The following statement is standard, see, e.g., \cite{Reu}, Theorems 3.1--3.2.

\begin{theorem}
\label{th:coproduct}

 Let $S\in\ov \U(\frg)$.

{\rm a)} $S\in \ov\frg$ if and only if $\delta(S)=S\otimes 1+1\otimes S$.

\smallskip

{\rm b)} $S\in \ov \G$ if and only if $\delta(S)=S\otimes S$.
\end{theorem}

 Identify $\ov \U(\frg)$ with the symmetric algebra
$S(\frg)$ in the usual way, see, e.g., \cite{Ser} or \cite{Dix}, Section 2.4.
Let $x^\alpha= x_1^{\alpha_1}\dots x_n^{\alpha_n}$ be a
 monomial. Denote
$p(x^\alpha)\in \U(\frg)$ the element obtained by symmetrization
of $x^\alpha$,
$$
p(x^\alpha)=
(c_\alpha)^{-1}\!\!\!\!\!\!\!\!\!\!\!
\sum_{\begin{matrix}\text{collections $\{j_1,j_2,\dots, j_{|\alpha|}\}$}
\\
\text{containing $\alpha_k$ entries of $k$}
\end{matrix}}
\!\!\!\!\!\!\!\!\!\!\!\!\!
 x_{j_1} x_{j_2}\dots x_{j_|\alpha|},
\qquad\text{where $c_\alpha=\frac{|\alpha|!}{\alpha_1!\dots\alpha_n!}$}
.
$$

 The following statement is obvious (see \cite{Dix}, 2.7.2).

\begin{lemma}
\begin{equation}
\frac 1{\alpha!}
\delta(p(x^\alpha))=\sum_{\beta,\gamma:\, \beta+\gamma=\alpha}
\frac 1{\beta!}p(x^\beta)\otimes\frac 1{\gamma!} p(x^\gamma)
\label{eq:shuffle}
\end{equation}
\end{lemma}

If $x_j$ is a basis in $\frg$, then elements $p_\alpha(x^\alpha)$
form a basis in $\U(\frg)$. Note that this produces a decomposition
$$
\ov\U(\frg)=\oplus_{j=0} \ov\U(\frg)_j
$$
according  degrees of polynomials $p$, Note that this decomposition differs from
the decomposition (\ref{z}) discussed above. Also the latter decomposition is 
a decomposition of a linear space and not a gradation of the enveloping algebra. 

Elements
$p_\alpha(x^\beta)\otimes p_\alpha(x^\gamma)$ form a basis
in $\ov\U(\frg)\otimes \ov\U(\frg)$.

Applying the identity (\ref{eq:shuffle})
 to a formal series $\sum c_\alpha p(x^\alpha)$, we 
immediately get an expansion of
$\delta\bigl(\sum c_\alpha p(x^\alpha)\bigr)$
in the basis. Emphasis that similar terms can not appear. More precisely

\begin{lemma}
\label{l:similar}
Let 
$$
T=
\sum_{k,l\ge 0}
t_{k,l},\qquad\text{ where $t_{k,l}\in\ov \U(\frg)_k\otimes \ov\U(\frg)_l$}.
$$
has the form $\delta(R)$, where
$$
R=\sum_{m\ge 0} r_m,\qquad\text{ where $r_m\in\ov \U(\frg)_m$}
.$$
 Then any term $t_{k,l}$ uniquely determines $r_{l+l}$.
\end{lemma}

{\sc Proof of Theorem \ref{th:coproduct}.}
 a) Having a term $p(x^\alpha)$
in $S$
 of degree $\ge 2$ we get a term in $\delta(S)$ of degree (1,1). 

 b). 
Let $S=\sum b_\alpha p(x^\alpha)$ satisfy
$\delta(S)=S\otimes S$. 
The expression  $\delta(S)$ contains a subseries
$\sum b_\alpha p(x^\alpha) \otimes 1 $. The expression
$S\otimes S$ contains a subseries  $\sum b_\alpha p(x^\alpha)\otimes b_0$.
Therefore,  $b_0=1$.
 We write
 $$
S=1+\sum \sigma_k x_k+\Bigl\{\text{terms $p(x^\alpha)$
of higher degree}\Bigl\}
. $$
Then the term of $S\otimes S$ of degree $(1,1)$ is
$\sum \sigma_k \sigma_l x_k\otimes x_l$. This term
must appear from $\delta(S)$. The only chance
is
\begin{multline*}
S=1+\sum \sigma_k x_k+ \frac 1{2!}\sum \sigma_k \sigma_l x_k x_l+\dots
=\\=
1+\sum \sigma_k x_k+\frac 1{2!}(\sum \sigma_k x_k)^2+\dots=
1+\sum \sigma_k x_k+ \frac 1{2!}\sum \sigma_k \sigma_l p(x_k x_l)+\dots
.
\end{multline*}
Next, we examine the term of $S\otimes S$ of degree
$(2,1)$. It has the form
$$
\frac 1{2!}(\sum \sigma_k x_k)^2\otimes (\sum \sigma_k x_k)	
$$
By Lemma \ref{l:similar}, the term of $S$ of degree 3 is $\frac 1{3!}(\sum \sigma_k x_k)^3$
(since this term is admissible and admissible term is unique).
 Repeating this argument we come to
$$
S=\sum \frac 1{k!} (\sum \sigma_k x_k)^k
.
$$


{\bf\punct Inverse element.} Recall that any enveloping algebra
$\ov\U(\frg)$ has a canonical anti-automorphism $\sigma$ determined by
$$\sigma(x)=-x \quad\text{for $x\in\frg$},
\qquad \sigma(z)\sigma(u)=\sigma(uz).$$

\begin{lemma}
\label{l:inverse}
For $S\in \ov\G$, we have $S^{-1}=\sigma(S)$.
\end{lemma}

{\sc Proof.} Indeed $S=\exp(x)$ for some $x\in\frg$. Then
$\sigma(S)=\sigma(\exp(x))=\exp(-x)=S^{-1}$.\hfill $\square$

\smallskip


{\bf\punct Example. Free Lie groups $\ov\bfFr_n$.}
 Consider the group
$\ov\bfFr_n$ corresponding to the free Lie algebra
$\frfr_n=\frfr[\omega_1,\dots,\omega_n]$.
Let $S\in \ov\bfFr_n$, $S=\sum_w c_w w$, where $w$ ranges in the set
of words in the letters $\omega_1$, \dots, $\omega_n$.
Then the equation $\delta(S)=S\otimes S$ is equivalent to the system
of quadratic equations
$$
c_v c_w=\sum_{\text{$u$, where $u$ is a shuffle of $v$ and $w$}} c_u
.
$$
Recall that a {\it shuffle} of $v$ and $w$ is a word
equipped with coloring of letters  in black and white
such that removing black letters we
get the word $v$ and removing white letters we get the word $w$.

On other characterizations of $\ov\bfFr_n$, see
 Reutenauer's book  \cite{Reu}.
 For instance, coefficients $c_v$ in the
front of Lyndon words%
\footnote{We introduce the alphabetic order $\prec$ on the set of all  monomials  in letters $\omega_j$.
A monomial $w$ is a {\it Lyndon word} if for any decomposition $w=uv$ we have $w\prec v$.}
 $v$ form a coordinate system on $\ov\bfFr_n$.


\smallskip

{\bf\punct Ordered exponentials.} Let $t$ ranges in a segment
$[0,a]$.
Let
$$
\gamma(t)=\sum_{j=1}^\infty \gamma^{[j]}(t),
$$
be a map $[0,a]\to \ov\frg$.

\begin{proposition}
\label{pr:exp1}
{\rm a)}
If each $\gamma^{[j]}(t)$ is an integrable bounded function,
then the differential equation
$$
\frac d{dt} \E(t)= \E(t)\gamma(t) ,\qquad \E(0)=1
$$
has a unique solution {\rm(}'ordered exponential'{\rm)} $\E(t)\in\ov \U(\frg)$.
Moreover, $\E(t)\in\ov \G$.

\smallskip

{\rm b)} Let $\gamma_m(t)$ be a sequence.
If for each $j $ the sequence  $\gamma_m^{[j]}$
 converges to $\gamma^{[j]}$ in the $L_1$-sense, then
$\E_m(a)$ converges to $\E(a)$.

\end{proposition}

{\sc Proof.} Expanding $\cE(t)=\sum \cE(t)^{[j]}$ we come to the system of equations,
$$
\Biggl\{ \frac d{dt}упорядоченная экспонента\E^{[j]}(t)=\gamma^{[j]}(t)+ \E^{[1]}(t)\gamma^{[j-1]}(t)
+\dots+ \E^{[j-1]}(t)\gamma^{[1]}(t)
$$
We consequently find $\E^{[1]}$, $\E^{[2]}$, \dots.
One can easily write  closed expressions for $\E^{[j]}(t)$
in terms of repeated integrals, such expressions imply b).

If $\mu(t)$ is piece-wise constant, then
its ordered exponent has the form
$$
\prod_{j=1} \exp(x_j)\in\ov \G
.$$
Next, we approximate $\gamma(t)$ by piece-wise function
and pass to a limit. By Theorem \ref{th:coproduct}, the group
$\ov\G$ is closed in $\ov\U(\frg)$.
\hfill $\square$

\smallskip


{\bf\punct Groups $\ov \bfG(\R)$ and $\ov\bfG(\C)$.} 
Consider the complexification $\frg_\C$ of the Lie algebra $\frg$.
We apply the same construction to the algebra $\frg_\C$ and get the group
$\ov \bfG(\C)$. Sometimes we also use the notation
$\ov \bfG(\R)$ for $\ov\bfG$.


\section{The Lie braid  group $\ov\bfBr_n$}

\COUNTERS

{\bf\punct A basis in $\U(\frbr_n)$.}
Consider the Kohno Lie algebra $\frbr_n$ with generators $r_{ij}$.
We say that a word of $r_{ij}$ is {\it good} if it has a form
\begin{equation}
w_1 w_2\dots w_{n-1}
\label{eq:word}
\end{equation}
where $w_1$ is a word composed of $r_{12}$, $r_{13}$, $r_{14}$ \dots, $r_{1n}$,
a word $w_2$ is  composed of $r_{23}$, $r_{24}$, \dots, $r_{2n}$,
etc. The word $w_{n-1}$ has the form $r_{(n-1)n}^k$.

\begin{proposition}
Good words form a basis in $\U(\frbr_m)$.
\end{proposition}

In other words we get a canonical bijection
$$
\U(\frbr_n)\leftrightarrow \Ass_{n-1}\otimes \Ass_{n-2}\otimes\dots \otimes \Ass_1
$$

{\sc Proof.}  By definition, $\U(\frbr_n)$ is the algebra with generators
$r_{ij}$ and quadratic relations
\begin{align}
r_{ij} r_{kl}&=r_{kl} r_{ij} \qquad \text{if $i$, $j$, $k$, $l$ are pairwise
distinct};
\label{eq:trans1}
\\
r_{ij}r_{jl}&=r_{jl} r_{ij}-r_{ij}r_{il}+r_{il} r_{ij}
\label{eq:trans2}
.
\end{align}
Let we have a word
$$
\dots r_{ij} r_{kl}\dots \in\U(\frbr_n)
$$
with $\min(i,j)>\min(k,l)$. Then we transform $r_{ij} r_{kl}$
according (\ref{eq:trans1}) or (\ref{eq:trans1}).
Removing all such disorders we come to a linear combination of good words.

Denote by $s_p$ the number of good words of degree $p$. Obviously,
$$
\sum s_p t^p=\prod_{j=1}^{n-1} (1-jt)^{-1}
.$$
Comparing with (\ref{eq:koh-poi}) we get that good words are linearly independent.
\hfill $\square$

\smallskip

The {\it rule of multiplication} is following.
Let we wish to multiply
$$
r_{ij}\cdot w_1 w_2\dots w_{n-1}
$$
(see (\ref{eq:word})), let $i<j$. We write
\begin{multline*}
r_{ij}\cdot w_1 w_2\dots w_{n-1}
=[r_{ij},w_1]w_2\dots w_{n-1} +\dots
\\
+
w_1 \dots [r_{ij}, w_{i-1}]w_i\dots w_{n-1}+
w_1 \dots  w_{i-1} (r_{ij} w_i)\dots w_{n-1}
\end{multline*}
and evaluate commutators
\begin{multline*}
[r_{ij},w_\alpha]=
[r_{ij}, w_{\alpha, k_1}w_{\alpha, k_2}\dots w_{\alpha, k_p}]
=\\=
[r_{ij},w_{\alpha, k_1}]w_{\alpha, k_2}\dots w_{\alpha, k_p}+
\dots+ w_{\alpha, k_1}w_{\alpha, k_2}\dots[r_{ij}, w_{\alpha, k_p}]
.
\end{multline*}
A bracket can be non-zero only if  $k_\theta=j$ or $\alpha=j$.
In such  case we get a sum of two good monomials.

\begin{corollary}
Let elements $S\in\U(\frbr_n)$ be written as linear combinations of
good monomials. Fix $\alpha$. Then the substitution
$r_{ij}=0$ for all  $i\le \alpha$ and $j\le \alpha$ determines
a well-defined homomorphism
$\U(\frbr_n)\to \U(\frbr_{n-\alpha})$.
\end{corollary}


{\bf\punct Lie braid groups.}

\begin{theorem}
\label{th:TTT}
An element $S\in\ov\U(\frbr_n)$ is contained in
the group $\ov\bfBr_n$ if and only if $S$ can be represented as a product
\begin{equation}
S= T_1 T_2\dots T_{n-1 }
,
\label{eq:TTT}
\end {equation}
where $T_j\in \ov\Ass[r_{j(j+1)},\dots r_{jn}]$
 is an element of the free Lie  group $\ov\bfFr_{n-j}$.
\end{theorem}

{\sc Proof.}
We apply Theorem \ref{th:coproduct}.
Using the canonical basis we identify spaces
$$\Xi:\ov\U(\frbr_n) \to \ov\U(\frfr_{n-1}\oplus\dots\oplus \frfr_{1})
.$$
Moreover, this identification is compatible with the co-product
$$
\delta\circ \Xi=(\Xi\otimes\Xi)\circ \delta
.
$$
For $\ov\U(\frfr_{n-1}\oplus\dots\oplus \frfr_{1})$ the required statement is obvious.
\hfill $\square$

\smallskip

Thus $\ov\bfBr_n$ is a product of its subgroups (not a direct product)
\begin{equation}
\ov\bfBr_n\simeq\ov\bfFr_{n-1}\times \ov\bfFr_{n-2}\times\dots\times \ov\bfFr_{1}
.
\label{eq:times}
\end{equation}
All subgroups
$
\ov\bfFr_{n-1}\times\dots\times \ov\bfFr_{n-j}
$
are normal and we have homomorphisms
$\ov\bfBr_{n}\to\ov\bfBr_{n-j}$ for all $j$.

\smallskip

The corresponding statement (see \cite{Xi}) on the level of Lie algebras is
$$
\frbr_n\simeq=\frfr_{n-1}\oplus \dots \oplus \frfr_1
$$
(recall again, this is not an isomorphism of Lie algebras)

\smallskip


{\bf \punct Decomposition  of ordered exponentials.} Our next remark: evaluation
of ordered exponentials in $\ov\bfBr_n$ is reduced to
successive evaluation
of ordered exponentials in the free Lie
 groups $\bfFr_1$, \dots, $\bfFr_{n-1}$. More precisely:

\begin{proposition}
Let $\gamma(t)$ be a way in $\frbr_n$. Decompose it
as
$$
\gamma(t)=\gamma_1(t)+\dots +\gamma_{n-1}(t), \qquad
\text{where $\gamma_j(t)\in\frfr_j= \frfr[r_{(n-j)(n-j+1)},\dots, r_{(n-j)n}] $}
$$
Then the solution $\E(t)$ of the differential equation
$$
\E'(t)=\E(t)\gamma(t),\qquad \E(0)=1
$$
admits a representation
\begin{equation}
\E(t)= U_{n-1}(t)\dots U_1(t)
\label{eq:razl}
,
\end{equation}
where $U_j(t)\in\ov\bfFr_j= \ov\bfFr[r_{(n-j)(n-j+1)},\dots, r_{(n-j)n}]$,
\begin{align*}
&U_1'(t)=U_1(t)\gamma_1(t), &U_1(0)=1
\\
&U'_2(t)=U_2(t) \cdot \bigl[U_1(t)\gamma_2(t)U_1(t)^{-1}\bigr], &U_2(0)=1
\end{align*}
\dots\dots
\begin{multline*}
U'_{n-1}(t)=U_{n-1}(t) \cdot
 \bigl[U_{n-2}(t)\dots U_1(t)\gamma_{n-1}(t)U_1(t)^{-1}
\dots U_{n-2}(t)^{-1}\bigr],
\\ U_{n-1}(0)=1
\end{multline*}
and
$$
U_1(t)\gamma_2(t)U_1(t)^{-1}\in\frfr_2, \quad 
U_2(t)U_1(t)\gamma_3(t)U_1(t)^{-1}U_2(t)^{-1}\in\frfr_3,\,\dots
$$
\end{proposition}

{\sc Proof.} For definiteness, set $n=4$.
Let us look for a solution in the form
(\ref{eq:razl}, by (\ref{eq:times} it exists,
$$
U_3'U_2 U_1+U_3 U_2' U_1+U_3 U_2 U_1'=
U_3 U_2 U_1 (\gamma_3+\gamma_2+\gamma_1)
.
$$
We apply the homomorphism $\ov\bfBr_4\to \ov\bfBr_2$ to
both sides of the equality and
get $U_1'=U_1\gamma_1$. After a cancellation, we come to
$$
U_3'U_2+U_3 U_2'=U_3 U_2 ( U_1\gamma_3U_1^{-1}+ U_1\gamma_2U_1^{-1})
.
$$
Next, we apply homomorphism $\ov\bfBr_4\to \ov\bfBr_3$
 to both sides of equality and  get $U_2'=U_2\cdot U_1\gamma_2U_1^{-1}$. 
After a cancellation, we find $U_3$. \hfill$\square$


\smallskip

{\bf\punct Embedding of the braid  group to $\ov\bfBr_n(\C)$.
\label{ss:KZ}}
Denote by $\Omega$ the space $\C^n$ without diagonals $\xi_i=\xi_j$.
The fundamental group  $\pi_1(\Omega)$ is $\Br_n$.
Consider the following differential 1-form on $\Omega$ with values in
$\ov\U(\frbr_n)$:
$$
d\omega =\sum_{k<l} r_{kl}\, d\ln(\xi_k-\xi_l)
$$
This form determines a flat connection (see Kohno \cite{Koh2}, \cite{Koh3}).
We take a reference point $\xi\in \Omega$ and for any
loop $\mu(t)$ starting at $\xi$ we write  the ordered exponential
of
$$
d\omega(\mu)= \sum_{k<l} r_{kl}\, d\ln\bigl(\mu_k(t)-\mu_l(t)\bigr)
\in \U(\frbr_n)
.
$$
Thus we get a homomorphism
\begin{equation}
\pi_1(\Omega)\simeq\Br_n\to \ov\bfBr_n(\C)
\label{eq:monodromy}
\end{equation}

\smallskip

{\sc Remarks.} a) 
 To find explicitly this homomorphism we must solve a system 
of linear differential equations with non-constant coefficients. As far as I know
an explicit solution is unknown.

\small

b)  Drinfeld \cite{Dri} (see also \cite{Mar}) obtained an imitation of this construction
on level of formal series, this gives  an embedding of the group $\Br_n$ 
to algebra of formal series of $\hbar\, r_{kl}$ (and therefore to $\ov\U(\frbr)$).
The author does not know does this embedding coincides with Malcev embedding.

\small

c) The subgroup $\Br_n\subset \ov \bfBr_n(\C)$ defined by  (\ref{eq:monodromy})
is contained also in 'small' groups $\bfBr_n^\circ(\C)$, $\bfBr_n(\C)^!$ defined in the next two
sections.


\section{Groups $\G^\circ$}

\COUNTERS

{\bf\punct Norm on free associative algebra.}
Let define {\it norm} on the $\Ass[\omega]$.
If $z=\sum c_w w$ is  the expansion of $\Ass[\omega]$
in basis monomials, then
$$
\|z\|:=\sum |c_w|
.$$
Obviously,
\begin{equation}
\|z+u\|\le\|z\|+\|u\|,\qquad \|zu\|\le\|z\|\cdot\|u\|,
\qquad\|[z,u]\|\le 2 \|z\|\cdot\|u\|
\label{eq:norm}
\end{equation}


{\bf\punct Norm on $\U (\frg)$.} Let us represent
$\frg$ as quotient of free Lie algebra $\frfr[\omega]$ with respect to
an ideal $I=\oplus I^{[j]}$, by our assumptions
$\omega_j$ constitute a basis of $\frg^{[1]}$. Denote by
$J$ the left ideal $J\subset \U(\frfr[\omega])$ generated
by $I$.

\begin{lemma}
{\rm a)} $J$ is a two-side ideal.

\smallskip

{\rm b)} $\frfr[\omega]\cap J=I$.
\end{lemma}

This is obvious.

\smallskip

We define {\it norm} on $\U(\frg)$ as the norm on the quotient
space $\U(\frg)/J$. In other words, for
$z\in \U(\frg)$ we consider all representations
of $z$ as noncommutative
polynomials in generators $\omega_j$,
$$
z=\sum c_w w,\qquad \text{where $w$ are of the form
$\omega_{j_1} \omega_{j_2}\dots \omega_{j_p}$}
$$
and set
$$
\|z\|=\inf \sum |c_w|
.
$$
This norm satisfies the same relations (\ref{eq:norm}).

We  use this norm only for homogeneous elements
of $\U(\frg)$.

\smallskip


{\bf \punct Algebra $\U^\circ(\frg)$. Definition.}
An element of $z=\sum z^{[p]}$ is contained in
$\U^\circ(\frg)$ if for any $C$
\begin{equation}
\lfloor z  \rceil_C:=
\sup_p \bigl(\|z^{[p]}\| e^{Cp}\bigr)<\infty
.
\label{eq:seminorm}
\end{equation}
Thus $\U^\circ(\frg)$ is a Frechet space%
\footnote{A Frechet space is a complete locally convex space, whose topology
is metrizable
 by a  translation-invariant metric.}
 with respect to the family
of semi-norms $\lfloor\cdot\rceil_C$.

\begin{observation} For 
$z\in \ov\U(\frg)$ define the series
\begin{equation}
\Phi_z(\xi):=\sum \|z^{[p]}\| \xi^p
.\label{eq:Phi}
\end{equation}
Then $z\in \U^\circ(\frg)$ if and only if $\Phi_z(\xi)$
 is an entire function of
the variable $\xi\in\C$.
\end{observation}

\begin{proposition}
\label{l:Ucirc}
{\rm a)} If $z$, $u\in \U^\circ(\frg)$, then $zu\in \U^\circ(\frg)$.

\smallskip

{\rm b)} The algebra $\U^\circ(\frg)$ does not depend on a choice of a basis 
$\omega_j$
in $\frg^{[1]}$.

\smallskip

{\rm c)} If $z\in \U^\circ_+(\frg)$ and $f(\xi)$ is an entire function,
 then $f(z)\in \U^\circ(\frg)$.

\smallskip

{\rm d)} If $z\in \U^\circ_+(\frg)$, then $\exp(z)\in \U^\circ(\frg)$.
\end{proposition}

{\sc Proof.}
a) We have
$$
\|(zu)^{[p]}\|\le \sum_{q,r: q+r=p} \|z\|^{[q]}\|u\|^{[r]}
$$
Therefore Taylor coefficients of $\Phi_{zu}(\xi)$ are dominated by the 
Taylor coefficients of $\Phi_z(\xi)\Phi_u(\xi)$.  Therefore $\Phi_{zu}$
is an entire function.

\smallskip

b) Let $\omega'_1$, \dots $\omega'_n$ be another basis in $\frg^{[1]}$.
Denote by $A$ the transition matrix, by $\|\cdot\|'$ another norm.
 Then
$$
\|A^{-1}\|^{-p}\|z^{[p]}\| \le \|z^{[p]}\|'\le \|A\|^p \|z^{[p]}\|
.$$

c)
For $f(\xi)=\sum a_m\xi^m$, we set $F(\xi)=\sum |a_m|\xi^m$.
We have
\begin{multline*}
\|f(z)^{[\alpha]}\|=
\bigl\|\sum_m a_m \bigl((z^{[1]}+z^{[2]}+z^{[3]}+\dots)^m\bigr)^{[\alpha]}\bigr\|
\le 
\\ \le
\sum_{m=1}^\alpha
\biggl(m!\, |a_m|
\sum\limits_{\sum s_j=m, \sum js_j=\alpha}
\frac{\|z^{[1]}\|^{s_1}}{s_1!}\dots
\frac{\|z^{[m]}\|^{s_m}}{s_m!}\biggr)
\end{multline*}
The last expression is the Taylor coefficient at $\xi^\alpha$
of an entire function $F(\Phi_z(\xi))$. \hfill
$\square$

\smallskip

{\bf \punct The group $\G^\circ$.} We define
$$
\G^\circ:=\ov \G\cap \U^\circ (\frg),\qquad
\frg^\circ:=\ov \frg\cap \U^\circ(\frg)
.$$
By Proposition \ref{l:Ucirc} $\G^\circ$ is closed with respect to
multiplications, by Lemma \ref{l:inverse}, $\G^\circ$ contains inverse elements.

By Proposition \ref{l:Ucirc} we have a well-defined map $\exp:\frg^\circ\to \G^\circ$.

\begin{observation}
Generally, this map is not surjective.
\end{observation}

Indeed, let $\frg=\frfr[\omega_1,\omega_2]$. We have
 $\exp(\omega_1)\exp(\omega_2)\in \G^\circ$. Opening
brackets in
$$
\ln\Bigl[
\exp(\omega_1)\exp(\omega_2)\Bigr]=\ln\Bigl[
\Bigl(1+\omega_1+\sum_{k>1} \frac{\omega_1^k}{k!}\Bigr)
\Bigl(1+\omega_2+\sum_{l>1} \frac{\omega_2^l}{l!}\Bigr)\Bigr]
,$$
we find a summand
$
(-\frac 1{2k}(\omega_1\omega_2)^k)
$
in the formal series, it is a unique summand with
$(\omega_1\omega_2)^k$.
\hfill $\square$

\smallskip

{\bf \punct The group $\bfBr_n^\circ$.}

\begin{proposition}
\label{pr:Brcirc}
Represent $S\in\bfBr_n^\circ$ as a product $S=T_{n-1}\dots T_1$
of elements $T_j\in\ov\bfFr_j$, see {\rm(\ref{eq:TTT})}. Then 
all factors $T_j\in \bfBr_n^\circ$. 
\end{proposition}

{\sc Proof} by induction.  Denote by $\pi$ the canonical homomorphism
$\ov\U(\frbr_n)\to \ov\U(\frbr_{n-1})$, this is a quotient
of $\ov\U(\frbr_n)$ by the ideal generated 
$r_{12}$, \dots, $r_{1n}$. Evidentely, for any homogeneous element
$z$, we have
$
\|\pi(z)\|\le \|z\|
$.
Indeed, if $z=\sum c_w w$ is a linear combination of monomials, then
$\pi(z)$ is obtained by removing of  monomials containing letters
$r_{12}$, \dots, $r_{1n}$. This diminish $\sum |c_w|$.

Hence, for $S\in\bfBr^\circ_n$ we have $\pi(S)\in\bfBr^\circ_{n-1}$.
Therefore, $T_{n-1}=S\pi(S)^{-1}\in \bfBr_{n}^\circ$.
\hfill$\square$

\smallskip

Thus $\bfBr_n^\circ$ is a semidirect product of its free Lie subgroups.
Obviously, 
$$
\bfBr_n^\circ\supset \bfFr_{n-1}^\circ\times\dots\times \bfFr_1^\circ.
$$
However, the author does not know is it $\supset$ or $=$.

\smallskip

{\bf\punct Ordered exponentials.}
Now let $\gamma:[0,a]\to \frg^\circ$ be a measurable
function, let for each $C$
$$
\sup_{t\in[0,a], j\in \N}
 \|\gamma^{[j]}(t)\|e^{Cj}<\infty.
 .$$

\begin{proposition}
Under this condition the solution $\E$
of the equation
\begin{equation}
\frac d{dt} \E(t)= \E(t)\gamma(t),\qquad \E(0)=1
\label{eq:ordexp2}
\end{equation}
is  contained in $\G^\circ$.
\end{proposition}

{\sc Proof.} Denote by $V_C$ the completion of $ \U(\frg)$
with respect to the seminorm (\ref{eq:seminorm}).

\begin{lemma} For $\gamma\in\frg^\circ$,
the linear operator
$$
L z= z\gamma
$$
 is bounded in the Banach space $V_C$.
\end{lemma}

{\sc Proof of Lemma.} Take $A>C$. We have
$\sup \|\gamma^{[q]}\| e^{Aq}<\infty$.
Therefore
$$
 \|(z \gamma )^{[r]}\|\le \sum_{p,q:\,p+q=r}
\|\gamma^{[q]}\| \,\| z^{[p]} \|\le \mathrm{const}
\cdot \sum_{p,q:\,p+q=r} e^{-Aq}e^{-Cp}\le\mathrm{const}
\cdot e^{-Cp}
.$$

{\sc Proof of the proposition.}
The differential equation (\ref{eq:ordexp2}) is
equivalent to the integral equation
\begin{equation}
\E(t)=1+\int_0^t \E(\tau) \gamma(\tau) \,d\tau
.\label{eq:volterra}
\end{equation}

The usual arguments (see any text-book on functional analysis, e.g.,
\cite{Tri}) show that this equation has a unique solution
in the Banach space of continuous functions
$[0,a]\to V_C$.

Since this is valid for all $C$, we get
$\E(t)\in\U^\circ(\frg)$. By Proposition \ref{pr:exp1},
$\E(t)\in\G^\circ$.
\hfill $\square$

\smallskip


{\bf\punct Representations of $\G^\circ$.}

\begin{theorem}
\label{pr:fin}
{\rm a)} Let $\rho$ be a representation
of $\frg$  in a finite-dimensional linear space $W$.
Then  $\rho$ can be integrated to a representation
of   $\G^\circ(\C)$  in $W$.

\smallskip

{\rm b)} Let operators $\rho(\omega_j)$ be anti-selfadjoint.
Then the group $\G^\circ(\R)$
acts in $W$ by unitary operators.
\end{theorem}

{\sc Proof.} a) The enveloping algebra $\U^\circ(\frg)$ acts in $W$.

\smallskip

b) We refer to Lemma \ref{l:inverse},
$\rho(S)^*=\rho(\sigma(S))=\rho(S^{-1})$.
 \hfill$\square$

\smallskip

\begin{proposition}
\label{pr:highest}
Let $H$ be a non-compact real simple Lie group admitting
highest weight representations%
\footnote{The list of such groups is: $\SU(p,q)$, $\Sp(2n,\R)$, $\SOS(2n)$,
$\mathrm{O}(n,2)$ and two real forms of $\mathrm E_6$ and  $\mathrm E_7$.}.
 Let $\rho_1$, \dots, $\rho_n$
be unitary highest weight representations of $H$.
 Then the group $\bfBr_n^\circ(\R)$ acts in
$\rho_1\otimes\dots \otimes \rho_n$ by unitary operators.
\end{proposition}

{\sc Proof.} 
We apply the  Knizhnik--Zamolodchikov construction
(see above \ref{ss:KZ1})
 and get a representation of $\frbr_n$ in the tensor
product. Decompose our tensor product with respect
to the diagonal subgroup $H\subset H\times\dots \times H$,
$$
 \rho_1\otimes\dots \otimes \rho_n=
\bigoplus_{\mu} (V_\mu\otimes \mu)
,
$$
where $\mu_j$ ranges in the set of highest weight representations
of $H$ and
$V_\mu$ are finite-dimensional spaces with trivial action of $H$.
The Lie algebra $\frbr_n$ acts by
$H$-intertwining operators, i.e.
in the finite-dimensional spaces $V_\mu$. We  apply  Theorem
\ref{pr:fin}.
\hfill $\square$


\section{Groups $\G^!$}

\COUNTERS

\smallskip

{\bf\punct Algebras $\U^!(\frg)$.} Fix $C$. Denote by $\U_a^!(\frg)$
the set of all elements $z\in  \ov\U(\frg)$ satisfying
$$
[[z]]_a=
\sup_p \|z^{[p]}\| \cdot p!\, a^{-p} <\infty
.$$
Note that $\U_a^!(\frg)$ is a Banach space with respect to the norm
$[[\cdot]]_a$.

We define
$\U^!(\frg)$ as
$$
\U^!(\frg)= \bigcup_a \U_a^!(\frg)
.$$
Thus, $\U^!(\frg)$
is an inductive limit of Banach spaces.

\begin{proposition}
\label{l:!}
{\rm a)} $\U^!(\frg)$ is an algebra.

\smallskip

{\rm b)} Moreover if $z\in \U_a^!(\frg)$, $u\in \U_b^!(\frg)$, then
$zu\in \U^!_{a+b}$ and
$$ [[zu]]_{a+b}\le [[z]]_a \cdot [[u]]_b  $$

{\rm c)} $\U^!(\frg)$ does not depend on a choice of a basis in
$\frg^{[1]}$.
\end{proposition}

{\sc Proof.} b)
Let $z\in \U_a^!(\frg)$, $u\in \U_b^!(\frg)$, i.e.,
$$
\|z^{[p]}\|\le \lambda \frac {a^p}{p!},\qquad
\|u^{[q]}\|\le \mu \frac {a^q}{q!}
.
$$
Therefore
$$
\|(zu)^{[r]} \|\le\lambda\mu \sum_{p+q=r}
\frac{a^p b^q}{p!\,q!}
=\frac {\lambda\mu }{r!} (a+b)^r
.
$$
Thus $zu\in \U_{a+b}^!(\frg)$.

\smallskip

c) See proof of Proposition \ref{l:Ucirc}.
\hfill $\square$

\smallskip

{\sc Example.}
Functions $e^{\omega_1\omega_2}$ and $e^{[\omega_1,\omega_2]}$
are not contained in $\U^!(\frfr_2)$.
\hfill $\square$

\smallskip


{\bf\punct Groups $\G^!$.}
We define the group $\G^!$ as
$$
\G^!:=\ov\G\cap \U^!(\frg).
$$
For any $a$ we define the group $\G_a^!$ as
the subgroup of $\G^!$ generated by the subset
$\ov\G\cap \U_a^!(\frg)$.  Evidently,
if $a>b$, then $\G_a^!\supset \G_b^!$, and
$\G^!=\cup \G^!_a$.%
\footnote{The author does not knows, is it $\G^!_a=\G^!$?}

\begin{proposition}
Represent $S\in\bfBr_n^!$ as a product $S=T_{n-1}\dots T_1$
of elements $T_j\in\ov\bfFr_j$, see {\rm(\ref{eq:TTT})}. Then 
all factors $T_j\in \bfBr_n^!$. 
\end{proposition}

See proof of Proposition \ref{pr:Brcirc}.

\smallskip


{\bf\punct Ordered exponentials.} Let $\mu(t)$, $t\in [0,T]$ be a continuously differentiable
curve in $\frg^{[1]}$,
$$
\mu(t)=\sum \mu_j(t)\omega_j.
$$

\begin{proposition}
\label{pr:exp-last}
{\rm a)}
The solution of the differential equation
\begin{equation}
\E'(t)=  \E(t)\mu(t),\qquad \E(0)=1
\end{equation}
is contained in $\G^!$.

\smallskip

{\rm b)} Moreover, the solution  is contained in $\cap_a \G^!_a$.
\end{proposition}

{\sc Proof.}
Let $\mu_\epsilon(t)$ be a piecewise constant function such that
$\|\mu(t)-\mu_\epsilon(t)\|<\epsilon$, let $\mu_\epsilon =c_j$ on
$(t_{j-1},t_{j})$. Denote by $\E_\epsilon(t)$  the solution
of the corresponding differential equation.
 As we have seen above (Proposition \ref{pr:exp1})
 $\E_\epsilon(t)$
converges in $\ov\U(\frg)$ to $\E(t)$.
If $t_{\alpha}<t<t_{\alpha+1}$, then we have
$$
\E_\epsilon(t)=
\exp\bigl\{ (t_1-t_0)c_1\bigr\}\exp\bigl\{ (t_2-t_1)c_2\bigr\}\dots
\exp\bigl\{(t-t_\alpha)c_{\alpha+1}\bigr\}
$$

Keeping in mind Proposition \ref{l:!} we get that $L_\epsilon(t)$ is contained in the unit
ball of the space $\U^!_{\sum \|c_j\|(t_j-t_{j-1})}(\frg)$
and all $L_\epsilon(t)$ are contained in the unit ball of a certain space
$\U^!_r(\frg)$. A  limit as $\epsilon\to 0$
 is contained in the same ball.

\smallskip

b) We divide $[0,T]$ into small segments, 
solve equations 
$$
\cE_j'(t)=\cE(t)\cdot \mu(t),\qquad \cE_j(t_j)=1
.$$
 Then
$\cE(T)=\cE_N(T)\dots\cE_1(t_2) \cE_0(t_1)$.
\hfill $\square$


\section{Actions of groups $\G^!_a(\R)$ on manifolds}

\COUNTERS


Let $\cM\subset \C^N$ be a closed complex  submanifold
satisfying the condition: if $\xi\in\cM$, then the complex conjugate
point $\ov \xi$ is contained in $\cM$. Let the intersection
$\cM\cap \R^N$ be a smooth compact manifold.

Let the Lie algebra $\frg$ acts by holomorphic vector fields
on $\cM$, denote by $\Omega_j$   vector fields corresponding to the
generators $\omega_j$ of $\frg$. Assume that $\Omega_j$
are real on $M$

Under these conditions we show that a certain group $\G^!_a(\R)$
 acts on $M$ by analytic diffeomorphisms%
\footnote{See examples in Introduction.}.

\smallskip

{\bf\punct Construction of the action.}
To any element of $z\in\U(\frg)$ we assign the differential operator
$\rho(z)$
on $\cM$,
$$
\rho(\omega_{i_1}\omega_{i_2}\dots \omega_{i_k})=
\Omega_{i_1}\Omega_{i_2}\dots \Omega_{i_k}
$$
For any element of $z=\sum z^{[p]}\in\ov\U(\frg)$ we assign the formal series
$\rho(z)=\sum \rho( z^{[p]})$ of differential operators.


\smallskip

 Denote by $\cA(\cM)$
the algebra of entire functions on $\cM$ equipped with the topology of uniform convergence on compact subsets. More generally, for an open subset $V\subset \cM$
we denote by $\cA(V)$ the algebra of analytic functions on $V$.

\begin{theorem}
\label{th}
{\rm a)} For sufficiently small $a$ for any
 $z\in \U^!_a(\frg)$ for any $F\in \cA(\cM)$
the formal series $\rho(z)F$  uniformly converges in a neighborhood
of $M$ in $\cM$.

\smallskip

{\rm b)} Moreover, for $z\in \ov\G(\R)\cap \U^!_a(\frg)$,
$$
\rho(z) F(m)=F(\cQ_z(\xi))
$$
for a certain analytic diffeomorphism $\cQ_z$ of $M$.

\smallskip

{\rm c)} Represent $z\in \G^!_a(\R)$ as $z=z_1 z_2\dots z_m$, where
$z_j\in \ov\G(\R)\cap \U^!_a(\frg)$. Then the formula
$$
\cQ_z:=\cQ_{z_m}\dots \cQ_{z_1}.
$$
determines a well-defined homomorphism from $\G^!_a(\R)$
to the group of analytic diffeomorphisms of $M$.

\smallskip

{\rm d)} Let $\cE(t)\in \cE$ be an ordered exponential as in Proposition {\rm\ref{pr:exp-last}}.
Let $\cE(a)=1$. Consider a family of diffeomorphisms $\cR(t)$ such that
$$
\cR(t)^{-1}\frac d{dt}\cR(t)=\rho(\mu(t)), \qquad 
\text{$\cR(0)$ is identical}
$$
Then the diffeomorphism $\cR(a)$ is identical.
\end{theorem}

{\sc Remark.}  Thus we get an action of an infinite-dimensional group on the manifold.
However this group is 'small'. The author does not know answer to the following questions.

a) Does the braid group $\Br_n$ is contained in the  'small' Lie group $\bfBr^!(\R)$? 
Does the braid group
acts on the space of Klyachko polygons? Are diffemorphisms corresponding to Drinfeld associators
are well defined?

b) Is it possible to find a completion  $\wt\frfr$
 of the free algebra $\frfr[\omega_1,\dots,\omega_k]$ such that  any 
collection of vector field $\Omega_1$, \dots, $\Omega_k$ on a compact manifold
generate a representation of $\wt\frfr$? If the answer is 'yes', is it possible to describe
the corresponding Lie group?
\hfill $\square$

\smallskip


{\bf\punct Estimate of derivatives.}

\begin{lemma}
\label{l:estim}
There exist neighborhoods $U\supset V$ of $M$ in $\cM$
and  constants  $c$, $A$ such that for any $F\in\cA(\cM)$
for any collection $j_1$, \dots, $j_p$
the following estimate holds
\begin{equation}
\sup_{m\in V}
\Omega_{j_1}\dots\Omega_{j_p} F(\xi)\le
A\cdot
p! c^p \sup_{\xi\in U} |F(\xi)|.
\label{eq:estim}
\end{equation}
\end{lemma}

{\sc Proof.}
Since the manifold $M$ is compact, the question is local.
So let us write the vector fields $\Omega_j$ in coordinates
$$
\Omega_j(\xi)=\sum_k^{\dim M} \alpha_{jk} \frac\partial{\partial \xi_k}
$$
We open brackets in $\Omega_{j_1}\dots\Omega_{j_p} F(\xi)$
and get $(\dim M)^p$ summands. It is sufficient to obtain an estimate of
the type (\ref{eq:estim}) for each summand. Hence we  assume
that each $\Omega_j$ has the form
$$
\Omega_j=\beta_j \frac\partial{\partial \xi_{k_j}}
.$$
Next, we define
operators
$$
H_j(t)= F(\xi_1,\xi_2,\dots)=F(\xi_1, \dots, \xi_{k_j}+t, \xi_{k_j+1},\dots)
$$
and write
$$
\Omega_{j_1}\dots\Omega_{j_p} F(\xi)=
\Bigl(\frac\partial{\partial t_1}\dots \frac\partial{\partial t_p}
H_{1}(t_1)\dots H_p(t_p) F(\xi)\Bigr)\biggr|_{t_1=\dots=t_p=0}
$$
Let  $F$, $\beta_1$, $\beta_2$, \dots
  be holomorphic in the polydisk $|\xi_\alpha|<\rho+\epsilon$. Then
$$g(\xi,t):=H_1(t_1)\dots H_p(t_p) F(\xi)$$
 is holomorphic
in the polydisk
$$
|\xi_1|<\rho/2, \quad |\xi_1|<\rho/2,\dots, |t_1|<\rho/2p,\dots,|t_p|<\rho/2p
.
$$
We estimate derivatives   by the multi-dimensional Cauchy formula
(see, e.g, \cite{GR}).
\begin{multline*}
\Bigl|\frac\partial{\partial t_1}\dots \frac\partial{\partial t_p}g(\xi,t)\bigr|_{t_1=\dots=t_p=0}\Bigr|
=\\=
\Bigl|\frac {(-1)^k}{(2\pi i)^k}
\int_{|t_1|=\rho/2p} \dots \int_{|t_1|=\rho/2p}
\frac{g(\xi,s)\,ds_1\dots ds_p}
{s_1^2\dots s_p^2}
\Bigr|\le C\cdot \Bigl(\frac {2p}{\rho}\Bigr)^p
\end{multline*}
under the assumption $|\xi_1|<\rho/2$, \quad $|\xi_1|<\rho/2$,\dots.
This is a desired estimate, recall that
$p!\sim \sqrt{2\pi p} \left(\frac pe\right)^p$.
\hfill $\square$


{\bf\punct Proof of Theorem \ref{th}.b.}
The statement a) follows from Lemma \ref{l:estim}.
More precisely, there is a neighborhood $V\subset \cM$ and $b$
such that for any $z\in \U^!_b(\frg)$ the operator
$\rho(z)$ is a well-defined map
$$
\rho(z):\cA(\cM)\to \cA(V)
.$$
Moreover, for $u$, $v\in \U^!_{b/2}(\frg)$
we have
$$
\rho(u)\rho(v)=\rho(uv)
.$$

The last equality is reduced to a permutation of absolutely convergent series.

\begin{lemma}
For $S\in \G^!\cap \ov \U_{b/2}$, $F_1$, $F_2\in \cA(\cM)$ we have
$$
\rho(S)(F_1\cdot F_2)=\rho(S)F_1\cdot \rho(S)F_2
.$$
\end{lemma}

{\sc Proof.} 
 Let $\delta:\ov\U(\frg)\to \ov\U(\frg)\otimes \ov\U(\frg)$
 be the co-product in enveloping
algebra as above. Let $\delta(z)=\sum x_i\otimes y_i$.
By the Leibnitz rule,
$$
\rho(z) (F_1F_2)=\sum \rho(x_i) F_1\cdot \rho(y_i) F_2
.$$
But $S$ satisfies  the condition $\delta(S)=S\otimes S$.
\hfill $\square$

\smallskip

Let $S\in\G^!\cap\ov\U_{b/2}(\frg)$.
Fix a point $v\in V$. We define a map
$\theta_v:\cA(\cM)\to\C$, by
$$
\theta_v(F)=\rho(S)F(v)
.$$
Then
$$
\theta_v(F_1+F_2)=\theta_v(F_1)+\theta_v(F_2), \qquad
\theta_v(F_1 F_2)=\theta_v(F_1)\theta_v(F_2)
.$$
Thus we get a character $\cA(\cM)\to \C$.  Since the manifold
$\cM$ is Stein%
\footnote{See the preamble to this section, on Stein manifolds see,
e.g., \cite{GR}.},
we get 
$$
\theta_v(F)=F(w) \qquad\text{for some $w=:\cQ_S(v)\in \cM$}
.$$
Therefore we get a map $\cQ_S: V\to \cM$.
Note that $F\circ \cQ_S$ is holomorphic for holomorphic
$F$ and
$$
F\circ \cQ_S\circ \cQ_{S^{-1}}=F=F\circ \cQ_{S^{-1}}\circ \cQ_S
.$$
Therefore $\cQ_S$ is a holomorphic embedding.

Next, $\cQ_S$ sends real functions on $M$ to real functions.
Therefore $\cQ_S$ preserves $M$.
\hfill $\square$

Moreover,
$$
 \cQ_{T}\cQ_{S}=\cQ_{ST}
\qquad
\text{ for
$S$, $T\in \G^!(\R)\cap\ov \U_{b/2}(\frg)$}
.
$$

\smallskip

{\bf\punct Proof of Theorem \ref{th}.c.}
The  statement c) of Theorem \ref{th}
is a corollary of the following lemma and the statement d) is a corollary of c).

\begin{lemma}
Let $z_1$, \dots $z_m\in \G^!_b(\R)$. Let
$z_1\dots z_m=1$ in $\G^!(\R)$.
Then
$
\cQ_{z_1}\dots \cQ_{z_m}=1
$.
\end{lemma}

{\sc Proof of lemma.}
For a non-zero $\tau\in \C^*$ we define the automorphism
$
A_\tau:\ov\U(\frg)\to\ov\U(\frg)
$
by
$$
A_\tau(z)=A_\tau\bigl(\sum z^{[p]}\bigr)
=\sum \tau^p z^{[p]}
.
$$
Evidentely, $A_\tau$ determines the automorphism of $\G^!(\C)$.

It can be easily checked that for $|\tau|<1$ the map
$\rho(A_\tau z):\cA(\cM)\to \cA(V)$ depends on $\tau$
holomorphically. Hence, for $S\in \G^!(\R)$
and $\tau\in(0,1)$ the map $\cQ_{A_\tau S}$ depends on $\tau$
real analytically.

 Let $\tau\in (0,1)$ and $z_j\in\G^!_b(\R)$ be as above. Then
$$
(A_\tau z_1)\dots (A_\tau z_1)=1
$$
On the hand $A_\tau z_j\in\G^!_{\tau b}(\R)$. Therefore, for
$\tau<b/m$ we have
$$
\cQ_{A_\tau z_1}\dots \cQ_{A_\tau z_m}=1
.$$
But the left hand side depends analytically on $\tau$, therefore
the identity holds
for all $t\in(0,1)$.
\hfill $\square$

{\tt Math.Dept., University of Vienna,

 Nordbergstrasse, 15,
Vienna, Austria

\&

Institute for Theoretical and Experimental Physics,

Bolshaya Cheremushkinskaya, 25, Moscow 117259,
Russia

\&

Mech.Math. Dept., Moscow State University,
Vorob'evy Gory, Moscow


e-mail: neretin(at) mccme.ru

URL:www.mat.univie.ac.at/$\sim$neretin

wwwth.itep.ru/$\sim$neretin
}

\end{document}